\documentclass[12pt,leqno]{article}
\usepackage{amsmath,amssymb,latexsym}
\usepackage{amscd}
\usepackage{mathrsfs}
\usepackage{srcltx}
\usepackage{theorem}
\usepackage{amsfonts}
\usepackage[notref,notcite]{showkeys}

\newcommand{\R}{\mathbb{R}}        
\newcommand{\C}{\mathbb{C}}        
\newcommand{\N}{\mathbb{N}}        
\newcommand{\Rn}{\mathbb{R}^d}     

\newcommand{\CN}{C^\infty(\Rn)}

\newcommand{\Exp}{{\rm Exp}}
\newcommand{\Log}{{\rm Log}}

\newcommand{\Di}{\cD(\Rn)}

\newcommand{\di}{\cD(\R)}
\newcommand{\Dip}{\cD\,'(\Rn)}

\newcommand{\dip}{\cD\,'(\R)}
\newcommand{\NZ}{\Rn_*}
\newcommand{\Dop}{\cD\,'(\Om)}
\newcommand{\Do}{\cD(\Om)}

\newcommand{\Proof}{\textbf{Proof:} \ }
\newcommand{\qed}{\hspace*{\fill} $\Box $}
\newcommand{\To}{\longrightarrow}

\newcommand{\cD}{\mathscr{D}}

\newcommand{\cE}{\mathscr{E}}
\newcommand{\cS}{\mathscr{S}}

\newcommand{\bn}{{\bf n}}

\newcommand{\Hb}{H_\bullet}

\newcommand{\om}{\omega}

\newcommand{\ssubset}{\subset\subset}

\newtheorem{proposition}{Proposition}
\newtheorem{example}[proposition]{Examples}
\newtheorem{lemma}[proposition]{Lemma}

\newtheorem{theorem}[proposition]{Theorem}

\newtheorem{definition}{Definition}

\newcommand{\proj}{\mbox{\rm proj}}

\newcommand{\Om}{{\Omega}}

\newcommand{\supp}{\mathrm{supp}\, }

\newcommand{\vp}{\varphi}

\title{{\sc Euler partial differential equations and Schwartz distributions.
}}
\author{Dietmar Vogt}
\date{}

\begin{document}

\maketitle

\footnotetext{
    \hfil\break\indent \begin{minipage}[t]{14cm}{\em Key words and phrases:} Euler differential operators, Schwartz distributions, distributions of finite order, global solvability. \end{minipage}.
     \hfil\break\indent
{}}

\begin{abstract}
 Euler operators are partial differential operators of the form $P(\theta)$
where $P$ is a polynomial and $\theta_j = x_j \partial/\partial x_j$. They are surjective
on the space of temperate distributions on $\Rn$. We show that this is, in general, not
true for the space of Schwartz distributions on $\Rn$, $d\ge 3$, for $d=1$, however, it is true.
It is also true for the space of distributions of finite order on $\Rn$ and on certain open
sets $\Om\subset\Rn$, like the euclidian unit ball.

\end{abstract}


Euler partial differential operators are operators of the form $P(\theta)$ where $P$ is a polynomial in $d$ variables and $\theta_j =x_j \partial/\partial x_j$ the Euler derivative. They are partial differential operators with variable coefficients which are singular at the coordinate hyperplanes. In contrary to partial differential operators with constant coefficients they admit distributional zero solutions with compact support, located at the singular locus. On $\CN$ they are surjective onto the annihilator of these zero solutions, in particular they are not surjective but have closed range (see \cite{DL1}). On the space $\cS'(\Rn)$ however they are surjective (see \cite{Vtemp}). The same holds on $\dip$ (see \cite{Vtemp} or Theorem \ref{t3} below). They have dense range in $\Dip$.  This raised the question whether they would be surjective on $\Dip$ for all $d$. In section \ref{s1} we show that this is not true, in general, at least for $d\ge3$. We give examples of differential equation $P(\theta)S=T$ without global solution $S$ even for such a natural operator as $P(\theta)=\sum_{j=1}^d \theta_j^2$.
The reason in this case is not a natural obstruction as in the $C^\infty$-case´, but impossibility to get a solution of locally finite order for some $T$ of infinite, but locally of finite order. In Section \ref{s2} we show that on the space $\cD'_F(\Rn)$ of distributions of finite order every non-trivial Euler differential operator is surjective, that is, every equation $P(\theta)S=T$ admits a global solution. The same holds for the open unit ball in $\Rn$. More general, we give for open $\Om$ conditions for global solvability on $\cD'_F(\Om)$ in terms of $P(\theta)$-convexity and we give conditions, which imply $P(\theta)$- convexity for every non-trivial polynomial $P$.

\section{Preliminaries}
We use the following notation $\partial_j=\partial/\partial x_j$, $\theta_j = x_j\,\partial_j$. For a multiindex $\alpha\in\N_0^d$ we set $\partial^\alpha = \partial_1^{\alpha_1}..\partial_d^{\alpha_d}$, likewise for $\theta^\alpha$. For a polynomial $P(z)=\sum_\alpha c_\alpha z^\alpha$ we consider the  \it Euler operator \rm $P(\theta)=\sum_\alpha c_\alpha \theta^\alpha$ and also the operator $P(\partial)$, defined likewise.

$P(\theta)$ and $P(\partial)$ are connected in the following way.  For $x\in\Rn$ we set
$\Exp(x)=(\exp(x_1),..,\exp(x_d)).$
$\Exp$ is a diffeomorphism from $\Rn$ onto $Q:=]0,+\infty[^d$. Its inverse is $\Log(x)=(\log(x_1),..,\log(x_d))$. The map
$C_\Exp:f\To f\circ\Exp$
is a linear topological isomorphism from $C^\infty(Q)$ onto $C^\infty(\Rn)$. For $f\in C^\infty(Q)$ we have
$P(\partial)(f\circ\Exp)=(P(\theta)f)\circ\Exp$ that is $P(\partial)\circ C_\Exp =C_\Exp\circ P(\theta)$. In this way solvability properties of $P(\theta)$ on $C^\infty(Q)$ can be reduced to solvability properties of $P(\partial)$ on $C^\infty(\Rn)$. This has be done in $\cite{V1}$ and essentially used in \cite{Vtemp}. There it was shown that every non-trivial Euler operator is surjective on $\cS'(\Rn)$, the space of temperate distributions. This result will be used in Section \ref{s2}. We will also make use of the fact that for elliptic $P(\partial)$ distributional zero solutions of $P(\theta)$ on some open $\Om\subset Q$ are real analytic on $\Om$.

Throughout the paper we use standard notation of Functional Analysis, in  particular, of distribution theory, and of the theory of partial differential operators. For unexplained notation we refer to \cite{DK}, \cite{H1}, \cite{MV}, 

\section{Examples for non-solvability}\label{s1}

We assume $d\ge 2$ and on $\R^{d+1}$ we use the variables $(x,y)$, $x\in\Rn$, $y\in\R$. We consider Euler differential operators  on $\cD'(\R^{d+1})$ and study the solvability of equations $P(\theta)S=T$ where
$$P(\theta)= Q(\theta_1,\dots,\theta_d) + \theta_y^p$$
and
$$T=\sum_{n=0}^\infty f_n(x)\otimes \delta^{(n)}(y).$$
with $f_n\in \Dip$ and $\supp f_n\subset \{x\,:\,|x|\ge n\}$.

We assume that $Q$ and the $f_n$ are chosen in such a way that every solution $F_n$ of $(Q(\theta)+(-n-1)^p) F_n=f_n$ does not vanish on any open subset of the `quadrant' $Q$. We might choose $Q(\theta)=\sum_{j=1}^d \theta_j^2$ and $f_n= \delta_\bn$ where $\bn=(n,\dots,n)$.

If we have such $F_n$ a natural candidate would be $S=\sum_n F_n(x)\otimes \delta^{(n)}(y)$ because
\begin{equation}\label{eq3} P(\theta)(F_n(x)\otimes \delta^{(n)}(y))=((Q(\theta) + (-n-1)^p)F_n(x))\otimes \delta_y^{(n)}.
\end{equation}
However, the series is not locally of finite order, hence does not define a distribution. Anyhow this shows the heart of the problem.

We study the problem on a strip around $\Rn$. We fix a function $\chi\in\cD([-2,+2]$ with $\chi(y)=1$ for $y\in[-1,+1]$ and consider functions $\vp(x,y)\in\cD(\R^{d+1})$ of the form
\begin{equation}\label{eq7}\vp(x,y)=\sum_{n=0}^m \vp_n(x)\chi(y)y^n/n!=\psi(x,y)\chi(y)\end{equation}
with
\begin{equation}\label{eq8}\psi(x,y)=\sum_{n=0}^m \vp_n(x)y^n/n!.\end{equation}
We obtain
$$P(\theta^*)\vp(x,y)=\sum_{n=0}^m (Q(\theta^*)+(-n-1)^p)\vp_n(x)\chi(y)y^n/n!+ (L\psi)(x,y)$$
where $L$ has the form
$$(L\psi)(x,y)= \sum_{j=1}^p(L_j\psi))(x,y)\chi^{(j)}(y).$$
In particular $\supp (L\psi)(x,y)\subset \{(x,y)\,:\,1\le |y| \le 2\}$.

Let now $S$ be a solution of $P(\theta)S=T$. We define $S_n\in\Dip$ by
 $$S_n(\gamma) = S(\gamma(x)\chi(y)y^n/n!).$$
 for $\gamma\in \Di$.

 We obtain
\begin{eqnarray*}
(P(\theta)S)\vp &=& \sum_{n=0}^m S_n((Q(\theta^*)+(-n-1)^p)\vp_n)+S(L\psi)
\\ &=& \sum_{n=0}^m ((Q(\theta)+(-n-1)^p) S_n)\vp_n+S(L\psi)\\
&=& \sum_{n=0}^m (-1)^n f_n(\vp_n).
\end{eqnarray*}
We define $R_n\in\Dip$ by
$$R_n = (-1)^n f_n -(Q(\theta)+(-n-1)^p) S_n)$$
and obtain for $\psi$ like in (\ref{eq8})
\begin{equation}\label{eq5}
\big(\sum_{n=0}^\infty (-1)^n R_n\otimes \delta^{(n)}\big)\psi = R(\psi)
\end{equation}
where
$$R(\psi)=S(L\psi)= \sum_{j=1}^p S((L_j\psi))(x,y)\chi^{(j)}(y)).$$

Both sides of (\ref{eq5}) define distributions in $(\Di\otimes \cE([-2,+2])'$ which coincide on the dense subspace of functions $\psi$ as in (\ref{eq8}), hence they coincide. The left hand side has support in $\Rn\times\{0\}$, while the right hand side has support in $\{(x,y)\,:\,1\le |y|\le 2\}$. So both sides are zero and we have
\begin{equation}\label{eq6}(Q(\theta)+(-n-1)^p)S_n=(-1)^n f_n \end{equation}
for all $n$. Returning to $\vp$ of the form like in (\ref{eq7}) we obtain
\begin{eqnarray*} S\vp &=& \sum_{n=0}^m S_n(\vp_n)= \sum_{n=0}^m S_n(\vp^{(0,n)}(x,0)) = \sum_{n=0}^m (-1)^n (S_n\otimes \delta^{(n)})\vp.
\end{eqnarray*}
So on the dense linear subspace of $\chi(y)\cdot\cD(\R^{d+1})$ consisting of functions $\vp(x,y)$ of the form (\ref{eq7}) we have
$$S = \sum_{n=0}^\infty (-1)^n S_n\otimes \delta^{(n)}$$
which implies that the sum must be locally finite contradicting our assumptions on the solutions of equation (\ref{eq6}). \qed

We have shown:
\begin{proposition}\label{t1} $P(\theta)$ as above is not surjective on $\cD\,'(\R^{d+1})$, $d\ge 2$.
\end{proposition}

This leads to the first main result of this paper.

\begin{theorem}\label{t2} For $d\ge 3$ there are Euler operators which are not surjective on $\Dip$.
\end{theorem}

We have shown non-surjectivity in the following cases.

\begin{example} { \rm For even $p$ and $d\ge 3$ the operator $\sum_{j=1}^d \theta_j^p$ is not surjective in $\Dip$. This holds, in particular for the ``Laplace-Euler''-operator $\sum_{j=1}^d \theta^2$. The same holds for $\sum_{j=1}^{d-1}\theta^2+i\theta_d$ and $\sum_{j=1}^{d-1}\theta^2+\theta_d$, the Euler-operators corresponding to the Schr\"odinger and to the heat equation and also for $\sum_{j=1}^{d-1}\theta_j^2-\theta_d^2$, the analogue to the wave equation. So for $d\ge 3$ we have counterexamples for the classical elliptic, parabolic and hyperbolic polynomials. }
\end{example}

For $d=1$ the situation is different (see the remarks at the end of \cite{Vtemp}).
\begin{theorem}\label{t3} Every non-trivial Euler-operator is surjective in $\dip$.
\end{theorem}

\Proof By the fundamental theorem of algebra it is enough to show it for $\theta-a$, $a\in \C$. For $T\in\dip$ we want to solve $(\theta-a)S=T$.

We set $\cD_0(\R)=\{\vp\in\di\,:\,\vp \text{ flat in } 0\}$. We find $S_0\in\cD_0(\R)'$ such that $x |x|^a S_0'=T_0 := T|_{\cD_0(\R)}$ and put $U_0=|x|^a S_0$. Then
$\theta U_0= a |x|^a S_0 + x |x|^a S_0' = a U_0+T_0$. Therefore $(\theta -a)U_0=T_0$. We extend $U_0$ by the Hahn-Banach theorem to $U\in\Dip$.
Then $\supp (\theta-a)U-T\subset \{0\}$. By \cite{Vtemp} we find $R\in\cS'(\R)$ such that $(\theta - a)R=(\theta-a)U-T$. Then $S=U-R$ solves the problem. \qed

\sc Remark: \rm Let us remark, that our counterexample does not work for $d=1$, that is, in $\cD\,'(\R^2)$. Because in this case the $F_n(x)$ (see the notation at the beginning of this section) can be chosen with support in half-lines beginning with $n$, hence the sum $S=\sum_n F_n(x)\otimes \delta^{(n)}(y)$ is locally finite and defines a distribution.

\section{Solvability in distributions of finite order}\label{s2}

In our examples of differential equations $P(\theta) S= T$ the distribution $T$ always was of infinite order and, of course, locally of finite order. We showed that there cannot exist a solution $S$ locally of finite order. In the proof in \cite{Vtemp}, that for temperate $T$ there is always a temperate solution $S$, one essential feature was that temperate distributions are always of finite order.  We use the result of \cite{Vtemp} to show that for $T$ of finite order there is always a solution $S$ of finite order. Moreover we develop a theory for arbitrary open subsets of $\Rn$. We should remark that the essential difficulty in \cite{Vtemp} was to handle the behaviour in the singular locus of the differential operator, that is, at the union of the coordinate hyperplanes. This was overcome in \cite{Vtemp} and we can use it here to provide local solvability.

From now on $P(\theta)$ is an arbitrary non-trivial Euler operator. $\Om$ and $\om$ denote open subsets of $\Rn$, $\cD_k'(\Om)$ the distributions of Sobolev order $k$ on $\Om$, see below.

We  follow the notation in \cite[\S 14]{MV}. We denote by $H^k$ the Sobolev space of order $k\in\N$. For open $\Om\subset\Rn$ the space $H_0^k(\Om)$ is the closure of $\Do$ in $H^k$. We set $\Hb^{-k}(\Om)=H^k_0(\Om)'$. Every $T\in \Hb^{-k}(\Om)$ can be extended to $\widetilde{T}\in (H^k)'\subset\cS'(\Rn)$. Due to \cite[Theorem 3.5]{Vtemp} and the open mapping theorem, applied to the surjective endomorphism $P(\theta)$ of $\cS'(\Om)$, we obtain:

\begin{lemma}\label{lem1} For every $k\in\N$ there is $m\in\N$ such that for every $\Om$ the following holds: for very $T\in\Hb^{-k}(\Om)$ there is $S\in\Hb^{-m}(\Rn)$ such that $P(\theta)S=T$ on $\Om$.
\end{lemma}

For every $k$ and $m=m_k$ chosen according to Lemma \ref{lem1} we set
$$E(\Om)=\{T\in \Hb^{-m}(\Om)\,:\,P(\theta)T\in \Hb^{-k}(\Om)\}.$$
We obtain an exact sequence
$$0\To N^m(\Om)\stackrel{j}{\To} E(\Om)\stackrel{P(\theta)}{\To} \Hb^{-k}(\Om)\To 0$$
where $N^m(\Om)=\{T\in\Hb^{-m}(\Om)\,:\,P(\theta)T=0\}$ and $j$ is the imbedding. This yields, due to reflexivity:
\begin{equation}\label{eq4}
N^m(\Om)'= E(\Om)'/P(\theta^*)H_0^k(\Om).
\end{equation}
We need some notation.
\begin{definition}\label{d1} For $\Om\subset\Rn$ open $\cD_F'(\Om)$ is the space of distributions of finite order. We set
$$\cD_k'(\Om)=\{T\in\Dop\,:\, T|_\om \subset \Hb^{-k}(\om)\text{ for all open }\om\ssubset\Om\}.$$
\end{definition}
Then we have $\cD_F'(\Om)=\bigcup_{k\in\N} \cD_k'(\Om)$.

\begin{definition}\label{d2}An open set $\Om\subset\Rn$ is called $P(\theta)$ convex if for every $k\in\N$ the following holds: for every $\om_1\ssubset\Om$ there is $\om_2\ssubset\Om$ such that for every $\om\ssubset\Om$ and $\vp\in H_0^k(\om)$ with $P(\theta^*)\vp\in H_0(\om_1)$ we have $\vp\in H_0(\om_2)$.
\end{definition}

Let now a $P(\theta)$-convex open set $\Om\ssubset\Rn$ be given. The we can find an exhaustion $\om_1\ssubset\om_2\ssubset..$ of $\Om$ such that for every $n\in\N$ the sets $\om_n\ssubset\om_{n+1}\ssubset\om_{n+2}$ are in the relation described in Definition \ref{d2}. We obtain a projective spectrum of exact sequences
$$0\To N^m(\om_n)\stackrel{j}{\To} E(\om_n)\stackrel{P(\theta)}{\To} \Hb^{-k}(\om_n)\To 0.$$
From equation (\ref{eq4}) it follows that every $\mu\in N^m(\om_n)'$ which vanishes on $N^m(\om_{n+2})$ also vanishes on $N^m(\om_{n+1})$,
and therefore $N^m(\om_{n+2})|_{\om_n}$ is dense in $N^m(\om_{n+1})|_{\om_n}$ in the topology of $N^m(\om_n)$. Then $P(\theta):\mathcal{E}(\Om)\to \cD_k'(\Om)$ is surjective, where $\mathcal{E}(\Om)=\{T\in\cD_m'(\Om)\,:\, P(\theta)T\in\cD'_k(\Om)\}=\proj_n E(\om_n)$ and, of course, $\cD_k'(\Om)=\proj_n \Hb^{-k}(\om_n)$.

The argument is standard. For the convenience of the reader we give the proof: for every $n$ we find $S_n\in E(\om_n)$ such that $P(\theta)S_n= T|_{\om_n}$. We set $R_1=R_2=0$ and determine inductively $R_n\in N^m(\om_n)$. Let $R_n$ be determined. Then $U_n=S_{n+1}-S_n+R_n \in N^m(\om_n)$ and for $n\ge 2$ we find $R_{n+1}\in N^m(\om_{n+1})$ such that $\|U_n-R_{n+1}\|_{\om_{n-1}}=\|(S_{n+1}-R_{n+1}) -(S_n-R_n)\|_{\om_{n-1}}\le 2^{-n}$. Clearly $S=\lim_n (S_n-R_n)$ exists everywhere on $\Om$ and $P(\theta)S=T$.

We have shown the second main result of this paper:
\begin{theorem}\label{t4} If $\Om\subset\Rn$ is open and $P(\theta)$-convex, then $P(\theta)$ is surjective in $\cD'_F(\Om)$.
\end{theorem}

To get examples of $P(\theta)$-convex sets we need some preparation. We set $Q=]0,+\infty[^d$. For $e\in\{-1,+1\}^d$ we set $Q_e=eQ$. We remark that for any $e$ we have $M_e\circ P(\theta)=P(\theta)\circ M_e$ where $M_e \vp(x)=\vp(ex)$. So the behaviour of $P(\theta)$ on $Q$ determines the behaviour on all `quadrants'.

A set $M\subset Q$ is called {\em m-convex} (cf. \cite{V1}) if $x^t y^{1-t}\in M$ for all $x,y\in M$ and $0<t<1$. $M\subset Q$ is m-convex if and only if $\Log \, M$ is convex. We call it strictly m-convex if $\Log \, M$ is strictly convex. A set $M\subset \NZ$ is called {\em strictly m-convex} if $e (M\cap Q_e)\subset Q$ is strictly m-convex for all $e\in\{-1,+1\}^d$.
 We obtain:

\begin{proposition}\label{p1} If $\Om$ has an exhaustion $\om_1\ssubset\om_2\ssubset\dots$ of open sets such that $\om_n\cap \NZ$ is strictly m-convex with $C^2$-boundary for all $n$, then $\Om$ is $P(\theta)$-convex for all $P(\theta)$.
\end{proposition}

\Proof On test functions $\vp$ we have $P(\theta^*) \vp=P(-1-\theta) \vp=:P^*(\theta)\vp$. Assume that $\vp\in H_0^k(\om_N)$ and $\supp P(\theta^*)\vp\subset \om_n$. Then $\tilde{\vp}:=\vp\circ \Exp$ is a function on $\Log \, (Q\cap \om_N)$ and $\supp P^*(\partial)\tilde{\vp} \subset \Log (Q\cap \om_n)$. By assumption $\Log \, (Q\cap \om_n)$ is strictly convex with $C^2$-boundary. Therefore it is, for any non-trivial $P$, intersection of non-characteristic half-spaces. By Holmgren's theorem (see \cite[Theorem 8.6.8]{H1}), $\supp \tilde{\vp} \subset \Log \,(Q\cap \om_n)$, hence $Q\cap \supp \vp \subset \om_n$. Applying this to $M_e \vp$ for all $e$ we obtain $\supp \vp\cap\NZ\subset \om_n$, hence $\supp \vp\subset \overline{\om}_n\subset \om_{n+1}$. \qed

\sc Remark: \rm The property of the $\om_n$ we really used was, that $\Log (e\,\om_n\cap Q)$ is an intersection of non-characteristic half spaces for all $n$.

Due to the concavity of $\log$ we get:

\begin{lemma}\label{lem2} If $M\subset Q$ is strictly convex and with any $y\in M$ and $x\le y$ (that is $x_j\le y_j$ for all $j$) also $x\in M$, then $M$ is strictly m-convex.
\end{lemma}

From Lemma \ref{lem2}, Proposition \ref{p1}, Theorem \ref{t4} and (for $p=\infty$) the Remark above we get the following examples.

\begin{example} For every non-trivial $P$ the Euler differential operator $P(\theta)$ is surjective on $\cD'_F(\Rn)$ and on $\cD'_F(\Om)$ where $\Om$ is the open unit ball of $\ell_p^d$, ($1\le p \le +\infty$), in particular, the open euclidian unit ball in $\Rn$.
\end{example}

\vspace{.5cm}

\noindent Bergische Universit\"{a}t Wuppertal,
\newline Fakult\"at f\"ur Mathematik und Naturwissenschaften,
\newline Gau\ss -Str. 20, D-42119 Wuppertal, Germany
\newline e-mail: dvogt@math.uni-wuppertal.de


\begin{thebibliography}{99}






\bibitem{DL1} P. Doma\'nski, M. Langenbruch, Surjectivity of Euler type differential operators on spaces of smooth functions, {\em to appear in Trans. Amer. Math. Soc.}

 \bibitem{DK} J. J. Duistermaat, J. A. C. Kolk, {\em Distributions. Theory and Applications}, Birkh\"auser, Boston, 2010.

 \bibitem{H1} L. H\"ormander, {\em The Analysis of Linear Partial Differential
Operators I}, Springer, Berlin-Heidelberg-New York-Tokyo, 1983.

 \bibitem{MV} R. Meise, D. Vogt: \it Introduction to functional analysis, \rm Clarendon Press, Oxford, 1997.



\bibitem{Seel} R. T. Seeley, Extension of $C^\infty$ functions defined in a half space, {\em Proc. Amer. Math. Soc.}, {bf 15} (1964), 615--626.


\bibitem{V1} D. Vogt, $\cE'$ as an algebra by multiplicative convolution,  {\em Functiones et Approximatio, DOI: 10.7169/facm/1719} (2018).


\bibitem{Vtemp} D. Vogt, Surjectivity of Euler operators on temperate distributions, {\em arXiv:1711.04140v2} (2018).


\end{thebibliography}
\end{document}